\documentclass{amsart} 
\usepackage{latexsym} 
\usepackage{amssymb}\usepackage{amsmath} 

\begin{document}


\newtheorem{theorem}{Theorem} 
\newtheorem{problem}{Problem} 
\newtheorem{definition}{Definition} 
\newtheorem{lemma}{Lemma} 
\newtheorem{proposition}{Proposition} 
\newtheorem{corollary}{Corollary} 
\newtheorem{example}{Example} 
\newtheorem{conjecture}{Conjecture} 
\newtheorem{algorithm}{Algorithm} 
\newtheorem{exercise}{Exercise} 
\newtheorem{remarkk}{Remark} 
 
\newcommand{\be}{\begin{equation}} 
\newcommand{\ee}{\end{equation}} 
\newcommand{\bea}{\begin{eqnarray}} 
\newcommand{\eea}{\end{eqnarray}} 
\newcommand{\beq}[1]{\begin{equation}\label{#1}} 
\newcommand{\eeq}{\end{equation}} 
\newcommand{\beqn}[1]{\begin{eqnarray}\label{#1}} 
\newcommand{\eeqn}{\end{eqnarray}} 
\newcommand{\beaa}{\begin{eqnarray*}} 
\newcommand{\eeaa}{\end{eqnarray*}} 
\newcommand{\req}[1]{(\ref{#1})} 
 
\newcommand{\lip}{\langle} 
\newcommand{\rip}{\rangle} 

\newcommand{\uu}{\underline} 
\newcommand{\oo}{\overline} 
\newcommand{\La}{\Lambda} 
\newcommand{\la}{\lambda} 
\newcommand{\eps}{\varepsilon} 
\newcommand{\om}{\omega} 
\newcommand{\Om}{\Omega} 
\newcommand{\ga}{\gamma} 
\newcommand{\rrr}{{\Bigr)}} 
\newcommand{\qqq}{{\Bigl\|}} 
 
\newcommand{\dint}{\displaystyle\int} 
\newcommand{\dsum}{\displaystyle\sum} 
\newcommand{\dfr}{\displaystyle\frac} 
\newcommand{\bige}{\mbox{\Large\it e}} 
\newcommand{\integers}{{\Bbb Z}} 
\newcommand{\rationals}{{\Bbb Q}} 
\newcommand{\reals}{{\rm I\!R}} 
\newcommand{\realsd}{\reals^d} 
\newcommand{\realsn}{\reals^n} 
\newcommand{\NN}{{\rm I\!N}} 
\newcommand{\DD}{{\rm I\!D}} 
\newcommand{\M}{{\rm I\!M}} 
\newcommand{\degree}{{\scriptscriptstyle \circ }} 
\newcommand{\dfn}{\stackrel{\triangle}{=}} 
\def\complex{\mathop{\raise .45ex\hbox{${\bf\scriptstyle{|}}$} 
     \kern -0.40em {\rm \textstyle{C}}}\nolimits} 
\def\hilbert{\mathop{\raise .21ex\hbox{$\bigcirc$}}\kern -1.005em {\rm\textstyle{H}}} 
\newcommand{\RAISE}{{\:\raisebox{.6ex}{$\scriptstyle{>}$}\raisebox{-.3ex} 
           {$\scriptstyle{\!\!\!\!\!<}\:$}}} 
 
\newcommand{\hh}{{\:\raisebox{1.8ex}{$\scriptstyle{\degree}$}\raisebox{.0ex} 
           {$\textstyle{\!\!\!\! H}$}}} 

\newcommand{\OO}{\won} 
\newcommand{\calA}{{\mathcal A}} 
\newcommand{\calB}{{\mathcal B}} 
\newcommand{\calC}{{\cal C}} 
\newcommand{\calD}{{\cal D}} 
\newcommand{\calE}{{\cal E}} 
\newcommand{\calF}{{\mathcal F}} 
\newcommand{\calG}{{\cal G}} 
\newcommand{\calH}{{\cal H}} 
\newcommand{\calK}{{\cal K}} 
\newcommand{\calL}{{\mathcal L}} 
\newcommand{\calM}{{\cal M}} 
\newcommand{\calO}{{\cal O}} 
\newcommand{\calP}{{\cal P}} 
\newcommand{\calS}{{\mathcal S}} 
\newcommand{\calU}{{\mathcal U}} 
\newcommand{\calX}{{\cal X}} 
\newcommand{\calXX}{{\cal X\mbox{\raisebox{.3ex}{$\!\!\!\!\!-$}}}} 
\newcommand{\calXXX}{{\cal X\!\!\!\!\!-}} 
\newcommand{\gi}{{\raisebox{.0ex}{$\scriptscriptstyle{\cal X}$} 
\raisebox{.1ex} {$\scriptstyle{\!\!\!\!-}\:$}}} 
\newcommand{\hna}{\hat{\nabla}}
\newcommand{\hpa}{\hat{\partial}}
\newcommand{\intsim}{\int_0^1\!\!\!\!\!\!\!\!\!\sim} 
\newcommand{\intsimt}{\int_0^t\!\!\!\!\!\!\!\!\!\sim} 
\newcommand{\pp}{{\partial}} 
\newcommand{\al}{{\alpha}} 
\newcommand{\sB}{{\cal B}} 
\newcommand{\sL}{{\cal L}} 
\newcommand{\sF}{{\cal F}} 
\newcommand{\sE}{{\cal E}} 
\newcommand{\sX}{{\cal X}}
\newcommand{\R}{{\rm I\!R}} 
\renewcommand{\L}{{\rm I\!L}} 
\newcommand{\vp}{\varphi} 
\newcommand{\N}{{\rm I\!N}} 
\def\ooo{\lip} 
\def\ccc{\rip} 
\newcommand{\ot}{\hat\otimes} 
\newcommand{\rP}{{\Bbb P}}
\newcommand{\rR}{{\Bbb R}} 
\newcommand{\bfcdot}{{\mbox{\boldmath$\cdot$}}} 
 
\renewcommand{\varrho}{{\ell}} 
\newcommand{\dett}{{\textstyle{\det_2}}} 
\newcommand{\sign}{{\mbox{\rm sign}}} 
\newcommand{\TE}{{\rm TE}} 
\newcommand{\TA}{{\rm TA}} 
\newcommand{\E}{{\rm E\,}} 
\newcommand{\won}{{\mbox{\bf 1}}} 
\newcommand{\Lebn}{{\rm Leb}_n} 
\newcommand{\Prob}{{\rm Prob\,}} 
\newcommand{\sinc}{{\rm sinc\,}} 
\newcommand{\ctg}{{\rm ctg\,}} 
\newcommand{\loc}{{\rm loc}} 
\newcommand{\trace}{{\,\,\rm trace\,\,}} 
\newcommand{\Dom}{{\rm Dom}} 
\newcommand{\ifff}{\mbox{\ if and only if\ }} 
\newcommand{\nproof}{\noindent {\bf Proof:\ }} 
\newcommand{\remark}{\noindent {\bf Remark:\ }} 
\newcommand{\remarks}{\noindent {\bf Remarks:\ }} 
\newcommand{\note}{\noindent {\bf Note:\ }}

\newcommand{\boldx}{{\bf x}} 
\newcommand{\boldX}{{\bf X}} 
\newcommand{\boldy}{{\bf y}} 
\newcommand{\boldR}{{\bf R}} 
\newcommand{\uux}{\uu{x}} 
\newcommand{\uuY}{\uu{Y}} 
 
\newcommand{\limn}{\lim_{n \rightarrow \infty}} 
\newcommand{\limN}{\lim_{N \rightarrow \infty}} 
\newcommand{\limr}{\lim_{r \rightarrow \infty}} 
\newcommand{\limd}{\lim_{\delta \rightarrow \infty}} 
\newcommand{\limM}{\lim_{M \rightarrow \infty}} 
\newcommand{\limsupn}{\limsup_{n \rightarrow \infty}} 
 
\newcommand{\ra}{ \rightarrow }

\newcommand{\ARROW}[1] 
  {\begin{array}[t]{c}  \longrightarrow \\[-0.2cm] \textstyle{#1} \end{array} } 
 
\newcommand{\AR} 
 {\begin{array}[t]{c} 
  \longrightarrow \\[-0.3cm] 
  \scriptstyle {n\rightarrow \infty} 
  \end{array}} 
 
\newcommand{\pile}[2] 
  {\left( \begin{array}{c}  {#1}\\[-0.2cm] {#2} \end{array} \right) } 
 
\newcommand{\floor}[1]{\left\lfloor #1 \right\rfloor} 
 
\newcommand{\mmbox}[1]{\mbox{\scriptsize{#1}}} 
 
\newcommand{\ffrac}[2] 
  {\left( \frac{#1}{#2} \right)} 
 
\newcommand{\one}{\frac{1}{n}\:} 
\newcommand{\half}{\frac{1}{2}\:} 
 
\def\le{\leq} 
\def\ge{\geq} 
\def\lt{<} 
\def\gt{>} 
 
\def\squarebox#1{\hbox to #1{\hfill\vbox to #1{\vfill}}} 
\newcommand{\nqed}{\hspace*{\fill} 
           \vbox{\hrule\hbox{\vrule\squarebox{.667em}\vrule}\hrule}\bigskip} 
 
\title{Malliavin Calculus for Degenerate Diffusions}

\author{ A. S. \"Ust\"unel} 
\maketitle 
\noindent 
{{\bf Abstract:}{\small{Let $(W,H,\mu)$ be the classical Wiener space
    on $\R^d$.  Assume that $X=(X_t(x))$ is a diffusion process satisfying
    the stochastic differential equation with diffusion and drift coefficients
  $\sigma: \R^n\to \R^n\otimes \R^d$, $b: \R^n\to \R^n$, $B$ is an $\R^d$-valued
  Brownian motion. We suppose that $b$ and $\sigma$ are Lipschitz. Let
  $P(x)$ be the orthogonal projection from $\R^d$ to its closed
  subspace $\sigma(x)^\star(\R^n)$, assuming that 
  $x\to P(x)$ is continuously differentiable,
  we construct a covariant  derivative $\hat{\nabla}$ on the paths of the diffusion
  process, along the elements of the Cameron-Martin space and
  prove that this derivative is closable on $L^p(\nu)$, where $\nu$
  represents the law of  the above diffusion process, i.e., $\nu=X(x)(\mu)$, the image of the
  Wiener measure under the function $w\to X_\cdot(w,x)$. We study the
  adjoint of this operator and  we prove several results:
  representation theorem for $L^2(\nu)$-functionals, the logarithmic
  Sobolev inequality
  for $\nu$. As applications of these results
the proof of the Logarithmic Sobolev inequality on the path space of
 Dyson's Brownian motion is given  by using the covariant derivative. We then
 explain how to use this theory for deriving the functional inequalities for the measures defined by the semigroups of the
 diffusion process at the time $t=1$ and with fixed starting point.  
 
.}}
\vspace{0.5cm}



\tableofcontents
\section{\bf{Introduction} }
\noindent
In a preceding work we have calculated a martingale representation
theorem for the functionals of degenerate It\^o  process  which is a unique weak solution of a stochastic differential equation with path depending coefficients
under the hypothesis of the  weak uniqueness, following the  path
breaking idea of C. Dellacherie, cf. \cite{CD}. In the 
calculations an important concept shows up, namely, the behaviour of
the random  projection operator which is defined as the orthogonal projection
from $\R^d$ to the closed subspace $\sigma^\star(t,X)(\R^n)$, where
$\sigma$ is the diffusion coefficient, this projection operator will be denoted as $P_s(X)$ if the diffusion has path dependent coefficients (non-Markov case) and by $P(X_s)$ in case the coefficients are Markov (i.e., the closed loop case, or the case where $b(t,X(w))=b(X_t(w) and \sigma(t,X(w))=
\sigma(X_t(w))$. In particular in the Markov process case we consider only the
homogeneous situation. This operator gives rise to the
It\^o representation theorem of the martingales adapted to the
filtration of the diffusion process. In fact we have shown that, for
any functional $F(X)$ defined on the paths of the diffusion process $X$,
which is square integrable w.r.t. to its law, there exists an $\R^d$-valued process
adapted to the filtration of $X$, say $\alpha(X)$, such that
$$
F(X)=E[F(X)]+\int_0^1 P(X_s)\alpha_s(X)\cdot dB_s\,,
$$
$\mu$-a.s. Although $\alpha$ is not unique, the process
$(P(X_s)\alpha_s(X),s\in [0,1])$ is unique $ds\times d\mu$-a.s. In the
case of strong solutions, we have also chaos representation \`a la
N. Wiener (cf.\cite{Ito}) for the functionals of $X$ using the basic
$\calF(X)$-martingale $m$ defined by 
$$
m_t=\int_0^tP_s(X)dB_s\,.
$$
To calculate the integrands in this chaos developement, as well as the
process $P(X)\alpha(X)$ above, we need to
create some kind of differential calculus on the paths of the
diffusion process $X$. The first kind of derivative which comes to
one's mind is some kind of extension of the $H$-derivative of Leonard
Gross (cf.\cite{LG_0,Kuo}), where $H$ denotes the Cameron-Martin space
associated to the
$\R^d$-valued Brownian motion which is at the origin of $X$, as it has
been done in the case of flat Malliavin Calculus,
cf. \cite{I-S}. However this approach fails since $H$-derivatives of
diffusion's functionals, even the simplest ones, fail to remain
measurable with respect to the sigma algebra generated by the
diffusion process, hence we can not express the representation results w.r.t.
the probability governing the process under investigation.
To circumvent this difficulty we propose a
{\bf{new derivation}} operator which consists of smoothing the
$H$-derivative w.r.t. the final instant sigma algebra, namely
$\calF_1(X)=\sigma(X_t,t\in [0,1])$, of the diffusion process. This
operator is denoted by $\hna$ with  a hat at the top to indicate the
role of the conditional expectation in its definition. It is then easy
to show that $\hna$, restricted to cylindrical diffusion functionals
$\calS(X)=\{f(X_{t_1},\ldots,X_{t_n}):t_1\leq t_2\leq\ldots \leq t_n\leq
1,\,f\in\calS(\R^n), n\geq 1\}$ is a derivation (under almost sure
equality) and it is a closable operator in $L^p(\nu)$ for $p>1$, where
$\nu$ denotes the law of $X$. We construct a path space functional analysis
based on this operator, its adjoint ( taken w.r.t. the law of
$X$)and the composition of it with the former, which is a kind of
Ornstein-Uhlenbeck operator with non-linear interactions.
Let us note that this is the first work where one has a nonlinear frame
(in the sense of Boson-Fermion situation) having a reasonable Wiener chaos
representation. We then  show several relations between  the
adjoint of this derivative operator , and the martingale representation
theorem; in particular we
calculate explicitly the integrand for  this representation theorem
which extends the flat case, which is
called in the flat case Clark-Haussman-Bismut-Ocone formula (cf. \cite{O1,O2}) . 
Let us note nevertheless that in the case of degenerate diffusions the
situation is quite different, for instance the derivative constructed
above is not even torsion free.(As an  immediate application of the
theory   in finance cf. \cite{CDU}). We prove here  the Poincar\'e inequality and the
logarithmic Sobolev inequality (associated to the derivative $\hna$)
with respect to $\nu$, which is the law of $X$.  We also  prove
these equalities for the path measure of a very singular process,
namely  the Dyson Brownian motion. We mean by degeneracy here the fact that the   probability law on 
the path space of this process  is singular with respect to the Wiener measure.
Afterwards we explain how to use all
the theory developed in this document to derive the Poincar\'e and
log-Sobolev inequalities for the semigroup measures of the diffusions
for fixed instants. As a final application, we construct  the ``normal''
process defined as the stochastic integral of the ortogonal of the projection
operator $P_s(X)$, namely $I_{\R^d}-P_s(X)$ and prove that conditionally w.r.t.
the final sigma algebra $\calF_1(X)$ of the diffusion process
$(X_t,t\in[0,1])$, it is a non-homogeneous Gaussian process with independent
increments. As an application of this fact, we prove a conditional integration
by parts formula which implies a 
conditional  Clark representation theorem and the associated  conditional logarithmic
Sobolev inequality for the paths of this normal process, which is new.

\section{\bf{Preliminaries and Notations }}
\noindent
$(W,H,\mu)$ represents  the classical Wiener space, i.e., $W=C([0,1],\R^d)$,
$H=\{\int_0^\cdot\dot{h}(s)ds: \dot{h}\in L^2([0,1],\R^d)\}$ is the Cameron-Martin space, with the norm defined by
$|h|_H=\|\dot{h}\|_{L^2([0,1],\R^d)}$ and $\mu$ is the Gauss measure on $W$, defined, for any $\alpha\in W^\star$ (dual of
$W$ under the topology of uniform convergence),
$$
\int_W\exp(\sqrt{-1}<\alpha,w>) \mu(dw)=\exp-\half|\tilde{\alpha}|_H^2\,,
$$
where $\tilde{\alpha}$ denotes the image of $\alpha$ under the injection $W^\star\to H$. We denote by $\nabla F$ the Sobolev derivative of the Wiener function $F$
defined as
$$
\nabla F=\sum_{i=1}^\infty \nabla_{e_i} F\,,
$$
where $(e_i,i\geq 1)$ is an orthonormal basis in $H$ 
as the $L^p(\mu)$-completion of the Gateaux derivative of cylindrical functions:
$$
\nabla_hF(w)=\lim_{\la\to 0}\frac{F(w+\la h)-F(w)}{\la}\,,
$$
for $h\in H$. We denote by $\DD_{p,k}$ the space of Wiener functions $F$ such
that
$$
\|F\|_{p,k}=\sum_{i=0}^k\|\nabla ^iF\|_{L^p(\mu,H^{\otimes i})}<\infty\,,
 $$
$\nabla^i$ denotes the $i$-th iterate of $\nabla$ and $\nabla^{\otimes i}F$
is regarded as an $H^{\otimes i}$-valued function. The operator $\delta$ is the formal adjoint of $\nabla$ w.r.t. the Wiener measure $\mu$ and it coincides with the It\^o integral on the vector fields $\xi: W\to H$ such that
$s\to\dot{\xi}(s)$ is adapted to the filtration of the Wiener space. We refer the reader to \cite{ASU}, \cite{ASU-1} and the references there for further
results about the linear case.

Let us  recall now  the representation theorem and its
consequences, which have all been proven in \cite{ASU-5}:
Let $X=(X_t,t\in [0,1])$ be a weak  solution of the following
stochastic differential equation:
\begin{equation}
\label{SDE}
dX_t=b(t,X)dt+\sigma(t,X)dB_t,\,X_0=x,
\end{equation}
where $B=(B_t,t\in [0,1])$ is an $\R^d$-valued  Brownian motion and
$\sigma: [0,1]\times C([0,1],\R^n)\to L(\R^d,\R^n)$ and $b:
[0,1]\times C([0,1],\R^n)\to \R^n$ are  measurable maps, adapted to the
natural filtration of $C([0,1],\R^n)$ and of linear
growth. Denote by $(\calF_t(X),t\in [0,1])$
 the filtration of $X$ and let us denote by $K$ the set of
$\R^n$-valued, $(\calF_t(X),t\in [0,1])$-adapted processes $\alpha(X)$,
s.t.
$$
E\int_0^1(a(s,X)\al_s(X),\al_s(X))ds<\infty\,,
$$
where $a(s,w)=\sigma(s,w)\sigma^\star(s,w),\,s\in[0,1],\,w\in C([0,1],\R^n)$. 
\begin{theorem}
\label{rep-thm}
The set $\Gamma=\{N\in
L^2(\calF_1(X)):\,N=E[N]+\int_0^1(\al_s(X),\sigma(s,X)dB_s),\,\al\in K\}$ is
dense in $L^2(\calF_1(X))$.
\end{theorem}
\noindent
Theorem \ref{rep-thm} permits us to show the following
{\bf{extension of the martingale representation
theorem to the functionals of degenerate diffusions}}:
\begin{theorem}
\label{adapt_thm}
Denote by $P_s(X)$ a measurable version of the orthogonal projection
from $\R^d$ onto $\sigma(s,X)^\star(\R^n)\subset\R^d$ and let $F\in
L^2(\calF_1(X))$ be any random variable with zero expectation. Then
there exists a process $\xi(X)\in L_a^2(dt\times dP;\R^n)$, adapted to
$(\calF_t(X),t\in [0,1])$, such that
$$
F(X)=\int_0^1 (P_s(X)\xi_s(X),dB_s)_{\R^d}=\int_0^1 (\xi_s(X),P_s(X)
dB_s)_{\R^n}
$$
a.s.\\
Conversely, any stochastic integral of the form
$$
\int_0^1(P_s(X)\xi_s(X),dB_s)_{\R^d}\,,
$$
where $\xi(X)$ is an $(\calF_t(X),t\in [0,1])$-adapted, measurable
process with $E\int_0^1|P_s(X)\xi_s(X)|_{\R^d}^2ds<\infty$, gives rise to an
$\calF_1(X)$-measurable random variable.
\end{theorem}

\remark Let $\eta$ be an adapted process such that $\eta_s$ belongs to
the orthogonal complement of $\sigma^\star(\R^n)$ in $\R^d$ $ds\times
dP$-a.s. Then $\eta+\xi$  can also be used to represent
$F(X)$. Hence $\xi(X)$ is not unique but $P(X)\xi(X)$ is always
unique.

As a consequence of above theorems, we get

\begin{theorem}
\label{iden-thm}
Let $\dot{u}\in L^2(dt\times dP,\R^d)$ be adapted to the Brownian filtration, then we have
\begin{equation}
\label{cond-exp}
E\left[\int_0^1(\dot{u}_s,dB_s)|\calF_1(X)\right]=\int_0^1(E[P_s(X)\dot{u}_s|\calF_s(X)],dB_s)
\end{equation}
almost surely.
\end{theorem}

\begin{remarkk}
If $K$ is a process adapted to the Brownian filtration with values in
$\R^m\otimes \R^d$ which is in $L^2(dt\times d\mu,\R^m\otimes \R^d)$,
then we obtain from Theorem \ref{iden-thm}
\begin{equation}
\label{mat_case}
E\left[\int_0^1 K_s
  dB_s|\calF_1(X)\right]=\int_0^1E[K_s|\calF_s(X)]P(X_s)dB_s
\end{equation}
$\mu$-almost surely.

\end{remarkk}

\begin{corollary}
\label{cor_1}
Let $h\in H^1([0,1],\R^d)$ (i.e., the Cameron-Martin space), denote by
$\rho(\delta h)$ the Wick exponential
$\exp(\int_0^1(\dot{h}_s,dB_s)-\half\int_0^1|\dot{h}_s|^2ds)$, 
 then we have
$$
E[\rho(\delta h)|\calF_1(X)]=
\exp\left(\int_0^1(P_s(X)\dot{h}_s,dB_s)-\half\int_0^1|P_s(X)\dot{h}_s|^2ds\right)\,.
$$
\end{corollary}

\begin{theorem}
\label{chaos-thm}
\begin{itemize}
\item Assume that $\calF_t(X)\subset \calF_t(B)$ for any $t\in [0,1]$, where
$(\calF_t(B),t\in [0,1])$ represents the filtration of the Brownian motion. Define the martingale $m=(m_t,t\in [0,1])$ as
$m_t=\int_0^tP_s(X) dB_s$,
then the set 
$$
K=\{\rho(\delta_m(h)):\,h\in H\}
$$ 
is total in
$L^2(\calF_1(X))$, where
$\rho(\delta_m(h))=\exp\left(\int_0^1(\dot{h}_s,dm_s)-\half\int_0^1|P_s(X)\dot{h}_s|^2ds\right)\,$.
In particular, any element $F$ of $L^2(\calF_1(X))$ can be written in
a unique way as the sum
\begin{equation}
\label{m_wi}
F=E[F]+\sum_{n=1}^\infty\int_{C_n}(f_n(s_1,\ldots,s_n),dm_{s_1}\otimes\ldots\otimes
dm_{s_n})
\end{equation}
where $C_n$ is the $n$-dimensional simplex in $[0,1]^n$ and $f_n\in
L^2(C_n, ds^{\otimes n})\otimes (\R^d)^{\otimes n}$.
\item   
More generally, without the hypothesis $\calF_t(X)\subset \calF_t(B)$,
for any $F\in L^2(\calF_1(X))\cap L^2(\calF_1(B))$, the conclusions of
the first part of the theorem hold true.
\end{itemize}
\end{theorem}

\begin{remarkk}
In this theorem we need to make the hypothesis $\calF_t(X)\subset
\calF_t(B)$ for any $t\in [0,1]$ to assure the non-symmetric chaos
representation (\ref{m_wi}). Without this hypothesis, although we have
Theorem \ref{adapt_thm} and when we iterate it we have a similar
representation, but, consisting of a finite number of terms. It is
not possible to push this procedure up to infinity since we have no
control at infinity.
\end{remarkk}
\begin{remarkk}
Let us note that if there is no strong solution to the equation
defining the process $X$, the chaotic  representation property may
fail. For example, let $U$ be a weak solution of 
\begin{equation}
\label{co-ex}
dU_t=\alpha_t(U)dt+dB_t\,,
\end{equation}
with $U_0$ given. Assume that (\ref{co-ex}) has no strong solution, as
it may happen in the famous example of Tsirelson (\cite{I-W}), i.e.,
$U$ is not measurable w.r.t. the sigma algebra generated by $B$, then
we have no chaotic representation property for the elements of
$L^2(\calF_1(U))$ in terms of the iterated stochastic integrals of
deterministic functions on $C_n, n\in \N$ w.r.t. $B$; the contrary
would imply the equality of $\calF_1(U)$ and of $\calF_1(B)$, which
would contradict the non-existence of strong solutions. 
\end{remarkk}

\section{\bf{Conditional and Covariant Sobolev derivatives  for degenerate
    diffusions' functionals}}
In this section, we shall work with the diffusion processes in the
time  homogeneous Markov case, i.e., $b(t,w,x)=b(x)$ and $\sigma(t,w,x)=\sigma(x)$, with
Lipschitz hypothesis and of linear growth.  To keep track of the nonlinearity due to
the diffusion's paths interactions, we need to develop a special
derivation w.r.t. the Cameron-Martin space of the defining Brownian
motion. The next result paves the road and it  is of fundemental
importance:
\begin{theorem}
\label{cond_thm}
For $h\in H$, define $\hna_h X_t=E[\nabla_hX_t|\calF_1(X)]$, then it
satisfies the following relation:
\begin{equation}
  \label{deriv_1}
\hna_h X_t=\int_0^t\partial\sigma (X_s) \hna_h X_s P(X_s)dB_s+\int_0^t
\partial b(X_s) \hna_h X_s ds+\int_0^t\sigma(X_s)\dot{h}(s)ds
\end{equation}
$\mu$-almost surely. If $\sigma$ and $b$ are twice differentiable and
bounded,  and if  the projection map is differentiable and of
linear growth, for the iterated derivative $\hna_h^2X_t$ we have the
following expression:
\beaa
\hna_h^2X_t&=&\int_0^t \partial^2\sigma(X_s) \hna_hX_s\otimes \hna_hX_s P(X_s)dB_s+\int_0^t\partial\sigma(X_s)\hna_h^2X_sP(X_s)dB_s\\
&&+\int_0^t \partial\sigma(X_s)\hna_hX_s \partial P(X_s)\hna_hX_s
P(X_s)dB_s+\int_0^t\partial \sigma(X_s)\hna_hX_s(P(X_s)+I_{\R^d})\dot{h}(s)ds\\
&&+\int_0^t\partial^2 b(X_s)(\hna_hX_s\otimes \hna_hX_s+\hna^2_hX_s)ds
\eeaa
$\mu$-almost surely.

\end{theorem}
\nproof
We have 
\begin{equation}
\label{deriv}
\nabla_h X_t=\int_0^t\partial \sigma (X_s) \nabla_h X_s dB_s+\int_0^t
\partial b(X_s) \nabla_h X_s ds+\int_0^t\sigma(X_s)\dot{h}(s)ds\,,
\end{equation}
taking the conditional expectation of both sides of this relation, the
term with stochastic integral follows from Theorem \ref{iden-thm} and
the Lebesgue integral part follows from the fact that 
$$
E[\nabla_hX_t|\calF_1(X)]=E[\nabla_hX_t|\calF_t(X)]
$$
due to the fact that $\calB_t=\sigma(B_s,s\leq t)$ is  independent of
the future increments of the Brownian motion (after $t$). For the second
derivative we just iterate the calculation of the first derivative.

\nqed

\begin{remark}
It is clear that $\hna_h^2X_t\neq E[\nabla_h^2X_t|\calF_1(X)]$.
\end{remark}

\begin{corollary}
Let $\calS(X)$ be the set of functions on $W$ defined as 
$$
\calS(X)=\left\{f(X_{t_1},\ldots,X_{t_m}): 0\leq t_1<\ldots< t_m,\,f\in
\calS(\R^{nm}),m\geq 1\right\}\,,
$$
where $\calS(\R^m)$ denotes the space of rapidly decreasing smooth
functions of Laurent Schwartz on $\R^m$. Let $h\in H$, assume that
$(F_k(X),k\geq 1)\subset \calS(X)$ converges to zero in
$L^2(\calF(X))$ and that $(\hna_hF_k(X),k\geq 1)$ is Cauchy in
$L^2(\calF(X))$, then 
$$
\lim_{k\to\infty}\hna_hF_k(X)=0
$$
$\mu$-a.s. In other words $\hna_h$ is a closable operator on
$L^2(X(\mu))$, where $X(\mu)$ denotes the law of the diffusion
process, or the image of the Wiener measure $\mu$ under the map $X$
which is defined by the diffusion process.
\end{corollary}
\nproof
Let $\eta=\eta(X)$ be the limit of $(\hna_hF_k(X),k\geq 1)$, then
using Theorem \ref{iden-thm}, we
have, for any cylindrical $G(X)\in \calS(X)$,
\beaa
E[\eta(X)\,G(X)]&=&\lim_kE[\hna_hF_k(X)G(X)]=\lim_kE[\nabla_hF_k(X)G(X)]\\
&=&\lim_kE[F_k(X)(-\nabla_hG(X)+G(X)\delta h)]\\
&=&\lim_kE[F_k(X)(-\hna_hG(X)+G(X) \delta(\rP(X) h))]=0\,,
\eeaa
where $\delta(\rP(X) h)=\int_0^1P(X_s)\dot{h}(s)dB_s$. As the cylinder functions of the form $G(X)$ are dense in
$L^p(\calF(X))$, $\eta(X)=0$ $\mu$-a.s., and this proves the closedness of $\hna_h$ hence that  of $\hna$.
\nqed

\noindent
For $F\in \calS(X)$, we define $\hna F$ as 
$$
\hna F=\sum_{i=1}^\infty \hna_{e_i}F\,,
$$
where $(e_i,i\geq 1)$ is an orthonormal basis in the Cameron-Martin
space $H$. The above theorem implies in particular that $\hna$ is a
closable operator from $L^p(\calF(X))$ into $L^p(\calF(X),H)$. We
shall denote its closure with the same notation. In particular, we
denote by $\M_{p,1}$ the completion of $\calS(X)$ w.r.t. the norm
$$
\|F(X)\|_{p,1}=\||\hna F(X)|_H\|_{L^p(\mu)}+\|F(X)\|_{L^p(\mu)}\,.
$$
For $p=2$ we use the following version:
$$
\|F(X)\|^2_{2,1}=E\left[|\hna F(X)|^2_H+|F(X)|^2\right]\,.
$$

\noindent
If $K$ is a separable Hilbert space, we denote by $\calS(X,K)$, the
set of $K$-valued cylindrical functions of the trajectories of $X$ by
replacing $\calS(\R^n)$ above by $\calS(\R^n)\otimes K$, where the
latter denotes the (projective) tensor product.

\begin{definition}
\label{div}
Let $\nu$ be the measure given as $X(\mu)$, i.e., the  law of the
diffusion process as a probability measure on the path space. We
denote the adjoint of $\hna$ w.r.t. $\nu$ as $\delta_\nu$ (or
sometimes as $\hna^\star$). Namely, for  $\xi(X)\in \calS(X)\otimes
H$, i.e., the $H$-valued cylindrical, $\calF(X)$-measurable functions
and $G(X)\in \calS(X)$, we write
\beaa
E[\delta_\nu\xi(X)\,G(X)]&=&E_\nu[\delta_\nu\xi\,G]\\
&=&E_\nu[(\xi,\hna G)_H]\\
&=&E[(\xi(X),\hna G(X))_H]\,.
\eeaa
\end{definition}
\begin{theorem}
\label{div_thm}
Let $\xi(X)\in L^2(\calF(X),H)$ be of the form
$\sum_{i<\infty}\xi_i(X)e_i$, where $\xi_i(X)\in\calS(X)$ and $e_i\in H$
for each $i\in\N$. 
\begin{enumerate}
\item
Then
$$
\delta_\nu\xi(X)=\sum_i\left[-\hna_{e_i}\xi_i(X)+\xi_i(X)\delta(\rP(X)e_i)\right]\,,
$$
where we define the action of $P(X)$ on the Cameron-Martin space $H$
as  
$$
(\rP(X)h)(t)=\int_0^tP(X_s)\dot{h}(s)ds\,,
$$
for $h\in H$.
In particular, we have
\begin{equation}
\label{c_exp}
E[\delta\xi (X)|\calF_1(X)]=\delta_\nu\xi(X)\,.
\end{equation}
\item If $\eta=\eta(X)\in \calS(X,H)$ is adapted to the filtration $(\calF_t(X),t\in
  [0,1])$,  where $\calS(X,H)$ denotes $H$-valued version of $\calS(X)$,  then
\be
\label{div_adap}
\delta_\nu\eta=\delta\rP(X)\eta
\ee
$\mu$-a.s.
\end{enumerate}
\end{theorem}
\nproof
Let $g(X)\in \calS(X)$, then, using Definition \ref{div} and the
integration by parts formula for $\mu$, we get
\beaa
E[(\xi(X),\hna g(X))_H]&=&E\left[\sum_i\xi_i(X)\hna_{e_i}g(X)\right]\\
&=&E\left[\sum_i\xi_i(X)\nabla_{e_i}g(X)\right]\\
&=&\sum_i E\left[-\nabla_{e_i}\xi_i(X)g(X)+g(X)\xi_i(X)\,\delta e_i\right]\\
&=&\sum_i E\left[-\nabla_{e_i}\xi_i(X)g(X)+g(X)\xi_i(X)\,E[\delta e_i|\calF_1(X)]\right]\\
&=&\sum_i E\left[-\hna_{e_i}\xi_i(X)g(X)+g(X)\xi_i(X)\,\delta(\rP(X) e_i)\right]\\
&=&E\left[\left(\sum_i \left(-\hna_{e_i}\xi_i(X)+\xi_i(X)\,\delta(\rP(X)
  e_i)\right)\right)g(X)\right]\,,
\eeaa
since $E[\delta e_i|\calF_1(X)]=\delta(\rP(X) e_i)$ as shown in Theorem \ref{iden-thm}.
The relation (\ref{c_exp}) holds then true for the cylindrical case,
the general case follows by passing to the limit in $\M_{2,1}(H)$
(i.e., the space $\M_{2,1}\otimes H$, completed Hilbert-Schmidt tensor
product).

To prove the second part, it follows already from the relation
(\ref{c_exp}) that $\delta \rP(X)\eta=E[\delta\eta|\calF_1(X)]$, let
now $G=G(X)\in \calS(X)$, then
\beaa
E[\delta_\nu\eta(X)\,G(X)]&=&E\left[(\eta(X),\hna G(X))_H\right]=
E\left[(\eta(X),\nabla G(X))_H\right]\\
&=&E[\delta\eta(X)\,G(X)]=E\left[E[\delta\eta(X)|\calF_1(X)]\,G(X)\right]\\
&=&E\left[\delta\rP(X)\eta\,G(X)\right]\,,
\eeaa
as, by Theorem \ref{iden-thm}, $\delta\rP(X)\eta$ is
$\calF_1(X)$-measurable, the proof follows.
\nqed

\begin{corollary}
\label{ibp_formula}
Let $F(X),G(X)\in \calS(X)$ and let $h\in H$, then we have the
integration by parts formula
\be
\label{ibp}
E[\hna_h F(X)\,G(X)=E\left[F(X)\,(-\hna_h G(X)+G(X)\delta \rP(X)h)\right]\,.
\ee
\end{corollary}
\nproof
As $E[\hna_h F(X)\,G(X)]=E[\nabla_hF(X)\,G(X)]$, 
from the Gaussian integration by parts formula, we have
\beaa
E[\hna_h F(X)\,G(X)]&=&E[F(X)(-\nabla_hG(X)+G(X)\delta h)]\\
&=&E\left[F(X)E\left[(-\nabla_hG(X)+G(X)\delta h)|\calF_1(X)\right]\right]\\
&=&E[F(X)(-\hna_h G(X)+G(X)\delta\rP(X) h)]\,.
\eeaa
\nqed

\begin{corollary}
\label{cor_div}
For $a(X)\in \M_{2,1}$ and $\xi(X)\in \calS(X,H)$,  we have
\be
\label{d_1}
\delta_\nu (a(X)\xi(X))=a(X)\delta_\nu\xi(X)-(\hna a(X),\xi(X))_H\,.
\ee
Moreover, if $\dot{\xi}(X)$ is adapted to the filtration
$(\calF_t(X),t\in [0,1])$, then
\be
\label{d_2}
\delta_\nu\xi(X)=\delta(\rP(X)\xi(X))\,.
\ee
In particular, for any $h,k\in H$, we have
\be
\label{d_3}
E[\delta(\rP(X)\nabla_h\rP(X)k)|\calF_1(X)]=\delta_\nu(\rP(X)\hna_h\rP(X)k)
\ee
$\mu$-a.s.
\end{corollary}
\nproof 
The only claim which is not immediate is the last one  and it follows
from Theorem \ref{iden-thm}.
\nqed

\begin{lemma}
\label{simple_lemma}
For any $h,k\in H$, we have
$$
\hna_k(\delta\rP(X)h)=(\rP(X)h,k)_H+\delta_\nu(\hna_k\rP(X)h)
$$
$\mu$-a.s.
\end{lemma}
\nproof
By definition
\beaa
\hna_k(\delta\rP(X)h)&=&E[\nabla_k\delta\rP(X)h|\calF_1(X)]\\
&=&(\rP(X)h,k)_H+E[\delta\nabla_k\rP(X)h|\calF_1(X)]\\
&=&(\rP(X)h,k)_H+\delta\rP(X)\hna_k\rP(X)h\\
&=&(\rP(X)h,k)_H+\delta_\nu\hna_k\rP(X)h\\
\eeaa
by Theorem \ref{iden-thm} and Corollary \ref{cor_div}.
\nqed

\begin{theorem}
\label{div2-calcul}
Let $\xi=\xi(X), \,\eta=\eta(X)$ be in $\calS(X,H)$, then the
following identity holds true:
\begin{eqnarray}
\label{div2-formula}
E[\delta_\nu\xi\,\delta_\nu\eta]&=&E[(\rP(X)\xi,\eta)_H+\trace(\hna\xi\cdot\hna\eta)]\\
&&+E[\trace(\hna\rP(X)\xi\cdot\hna\eta)+\trace(\hna\rP(X)\eta\cdot\hna\xi)]\nonumber
\end{eqnarray}
\end{theorem}
\nproof
Writing
$E[\delta_\nu\xi\,\delta_\nu\eta]=E[(\xi,\hna\delta_\nu\eta)_H]$, then
we shall calculate $\hna\delta_\nu\eta$. Let $k\in H$ and let $(e_i,i\geq
1)$ be an orthonormal basis of $H$, then
\beaa
\hna\delta_\nu\eta&=&\hna_k\left[\sum_i\eta_i\delta\rP(X)e_i-\hna_{e_i}\eta_i
\right]\\ 
&=&\sum_i\left[\hna_k\eta_i\delta\rP(X)e_i+\eta_i\hna_k\delta\rP(X)e_i-\hna_k\hna_{e_i}\eta_i\right]\\
&=&\sum_i\left[\hna_k\eta_i\delta\rP(X)e_i+\eta_i(\rP(X)e_i,k)_H\right]\\
&&+\sum_i\left[\eta_i\delta_\nu\hna_k\rP(X)e_i-\hna_k\hna_{e_i}\eta_i\right]
\eeaa
Taking the expectations after replacing the vector $k$ with $e_k$ and
using the above equality, we get

\beaa
E[\delta_\nu\xi\,\delta_\nu\eta]&=&E[(\xi,\hna\delta_\nu\eta)_H]
=E\left[\sum_i\eta_i\delta_\nu(\delta\rP(X)e_i)+(\rP(X)\xi,\eta)_H\right]\\
&&+E\left[\sum_{k,i}\eta_i(\delta_\nu\hna_{e_k}\rP(X)e_i)\xi_k-\xi_k\hna_{e_k}\hna_{e_i}\eta_i\right]\,,
\eeaa
where $\eta_i=(\eta,e_i)_H,\xi_i=(\xi,e_i)_H$. Using the relation
$$
\hna_\xi\hna_{e_i}\eta_i=\hna_{e_i}\hna_\xi\eta_i-\hna_{\hna_{e_i}\xi}\eta_i\,,
$$
we have

\beaa
E[\delta_\nu\xi\,\delta_\nu\eta]&=&E\left[\sum_i\hna_{\xi}\eta_i\delta\rP(X)e_i+(\rP(X)\xi,\eta)_H\right]\\
&&+E\left[\sum_{k,i}\eta_i\delta_\nu(\hna_{e_k}\rP(X)e_i)\xi_k-\hna_{e_i}\hna_\xi\eta_i+\trace(\hna\eta\cdot\hna\xi)\right]\\
&=&E\left[\delta_\nu\hna_\xi\eta+(\rP(X)\xi,\eta)_H+\trace(\hna\eta\cdot\hna\xi)+\sum_{k,i}\eta_i\delta_\nu(\hna_{e_k}\rP(X)e_i)\xi_k\right]\,.
\eeaa
Using the identity
$$
\delta_\nu(g\,\tau)=g\delta_\nu\tau-(\tau,\hna g)_H\,,
$$
for $g\in \calS(X),\,\tau\in \calS(X,H)$, we have
\bea
\label{div_stuff}
&&\sum_{k,i}E\left[\eta_i\delta_\nu(\hna_{e_k}\rP(X)e_i)\xi_k\right]\\
&=&\sum_{k,i}
E\left[\eta_i\delta_\nu(\xi_k\hna_{e_k}\rP(X)e_i)+\eta_i(\hna\xi_k,\hna_{e_k}\rP(X)e_i)_H\right]\nonumber\\
&=&E\left[\sum_i\eta_i\delta_\nu(\hna_\xi\rP(X)e_i)+\sum_k(\hna\xi_k,\hna_{e_k}\rP(X)\eta-\rP(X)\hna_{e_k}\eta)_H\right]\nonumber
\eea
and since
\be
\label{open}
\eta_i\delta_\nu(\hna_\xi\rP(X)e_i)=\delta_\nu(\eta_i\hna_\xi\rP(X)e_i)+(\hna_\xi\rP(X)e_i,\hna\eta_i)_H\,,
\ee
we obtain by substituting (\ref{open}) in (\ref{div_stuff})
$$
\sum_{k,i}E\left[\eta_i\delta_\nu(\hna_{e_k}\rP(X)e_i)\xi_k\right]=I+II\,,
$$
where
\beaa
I&=&\sum_i
E\left[\delta_\nu(\eta_i\hna_\xi\rP(X)e_i)+(\hna_\xi\rP(X)e_i,\hna\eta_i)_H\right]\\
&=&E[\trace(\hna_\xi\rP(X)\cdot\hna \eta)]
\eeaa
and
\beaa
II&=&\sum_{k,i}E[\eta_i(\hna\xi_k,\hna_{e_k}\rP(X)e_i)_H]\\
&=&\sum_kE[(\hna\xi_k,(\hna_{e_k}\rP(X))\eta)_H\\
&=&E[\trace(\hna\xi\cdot (\hna\rP(X))\eta)]
\eeaa
and this completes the proof.
\nqed

\noindent
The next result shows the local character of the derivative operator
$\hna$, which is well-known in the flat case:

\begin{proposition}
  \label{loc_prop}
  Assume that $F=F(X)$ is in $\M_{2,1}$, then $\hna F=0$ on the set
  $\{F=0\} $ $\nu$-almost surely.
\end{proposition}
\nproof
Let $\theta$ be a positive, smooth function of compact support on
$\R$ with $\theta(0)=1$,  denote by $\zeta$ its primitive and let
$\zeta_\eps(t)=\eps\zeta(t/\eps)$. Clearly $\zeta_\eps\circ F$ belongs
to $\M_{2,1}(X)$. If $u\in \calS(X,H)$, from the
definitions
\beaa
E[\zeta_\eps(F)\,\delta_\nu u]&=&E[(\hna\zeta_\eps(F),u)_H]\\
&=&E[\theta(F/\eps)(\hna F,u)_H]\to E[1_{\{F=0\}}]\,.
\eeaa
On the other hand, from the Dominated Convergence Theorem, we have
$\lim_{\eps\to 0}E[\zeta_\eps(F)\delta_\nu u]=0$ and as $\calS(X,H)$ is
dense in $\M_{2,1}(X,H)$, the proof follows.
\nqed

\noindent
The following result shows that $\hna$ defines a Markovian Dirichlet
form on the path space of the diffusion process:
\begin{proposition}
  \label{min}
  Assume that $f,\,g\in \M_{2,1}$, then $f\wedge g\in \M_{2,1}$, where
  $\wedge$ denotes the minimum.
\end{proposition}
\nproof
Let $(f_n,n\geq 1)$ and $(g_n,n\geq 1)$ be two sequences of
cylindrical functions approximating $f$ and $g$ in $\M_{2,1}$
respectively. We have
\beaa
\hna(f_n\wedge g_n)&=&E[\nabla f_n\, 1_{\{f_n\geq g_n\}}+\nabla g_n\,
1_{\{f_n<g_n\}}|\calF_1(X)]\\
&=&\hna f_n 1_{\{f_n\geq g_n\}}+\hna g_n 1_{\{f_n< g_n\}}\,,
\eeaa
Hence
$$
|\hna( f_n\wedge g_n)|_H\leq |\hna f_n|_H+|\hna g_n|_H\,.
$$
Consequently, $(f_n\wedge g_n,n\geq 1)$ is bounded in $\M_{2,1}$, then
it has a subsequence which converges weakly. This implies that
$f\wedge g\in \M_{2,1}$ and that $|\hna (f\wedge g)|_H\leq |\hna
  f|_H+|\hna g|_H$ almost surely.
  \nqed

\noindent
While making calculations we need also a second covariant derivation
on the paths of the diffusion process, this time this will be the
smoothed derivative w.r.t. the state space elements, namely we define
$$
\hpa_\xi X_t(\xi)=E[\partial_\xi X_t(\xi)|\calF_1(X(\xi))]
$$
where $\partial_\xi$ denotes the derivative of the diffusion process at
the instant $t\in [0,1]$ with respect to its initial starting point
$\xi\in \R^n$, similarly, for smooth functions on $\R^n$, we define
$$
\hpa_\xi f(X_t(\xi))=E[\partial f(X_t(\xi))\partial_\xi X_t(\xi)|\calF_1(X(\xi))]=
\partial f(X_t(\xi))\hpa_\xi X_t(\xi),
$$
in such a way that $\hpa$ behaves as a derivative w.r.t. state space variable of the functionals
of the paths of the diffusion process. We have the immediate result
which follows from Theorem \ref{iden-thm}

\begin{theorem}
  \label{J-calcul}
  The state space covariant derivative  $\hpa_xX_t(x)=J_t(x)$ satisfies
  the following  stochastic differential equation:
  \begin{equation}
    \label{J-eqn}
    J_t=I_{\R^n}+\int_o^t\partial\sigma(X_s(x))J_s(x)P(X_s(x))dB_s+\int_0^t\partial
    b(X_s(x))J_s(x)ds\,,
  \end{equation}
  $P$-a.s.,  for any $t\in [0,1]$.
\end{theorem}
\nproof
We have
\begin{equation}
  \label{der_1}
\partial X_t(x)=I_{\R^n}+\int_o^t\partial\sigma(X_s(x))\partial X_s(x)dB_s+\int_0^t\partial
b(X_s(x))\partial X_s(x)ds\,.
\end{equation}
Taking the conditional expectation  of both sides of the equality (\ref{der_1}),
we get
$$
E\left[\int_o^t\partial\sigma(X_s(x))\partial X_s(x)dB_s|\calF_1(X)\right]
=\int_o^tE[\partial\sigma(X_s(x))|\calF_s(X)]P(X_s(x))dB_s
$$
due to Theorem \ref{iden-thm}. We also have
$$
E\left[\int_0^t\partial b(X_s(x))\partial X_s(x)ds|\calF_1(X)\right]
=\int_0^t\partial
b(X_s(x))E[\partial X_s(x)|\calF_s(X)]ds
$$
from the Markov property of $(X_t(x))$ w.r.t. the Brownian filtration.
\nqed

\noindent
It follows from Theorem \ref{J-calcul} and from the It\^o formula that
\begin{proposition}
  \label{J-1-calcul}
  $(J_t,t\in [0,1])$ is a semimartingale with values in $GL(n,\R)$
  (non-singular matrices) and its inverse, denoted as $(K_t,t\in
  [0,1])$ has the following representation:
  \begin{eqnarray}
    \label{K-eqn}
    dK_t&=&-\sum_{i=1}^dK_t\partial\sigma^i(X_t(x))dm^i_t\\
     &&+K_t\left[\sum_{i,j=1}^d \partial\sigma^i(X_t(x))
        \partial\sigma^j(X_t(x))P_{i,j}(X_t(x))-\partial b(X_t(x))\right]dt\,,\nonumber
    \end{eqnarray}
where $dm_t=P(X_t(x))dB_t$, $m_t=(m^1_t,\ldots,m^d_t)$, and $\sigma^i$
is the $i$-th column vector of the $n\times d$-matrix $\sigma$.
\end{proposition}
\nproof
The equation (\ref{K-eqn}) is a linear stochastic differential
equation with bounded coefficients, hence it has a unique solution
without explosion. To prove the claim it suffices to show that
$J_tK_t=K_tJ_t=I_{\R^n}$ almost surely, where $I_{\R^n}$ denotes the identity map of
$\R^n$,  and this last equation follows from the It\^o formula.
\nqed

\noindent
An immediate consequence of these calculations, combined with the
variation of constants' method, is the following
\begin{corollary}
  \label{calcul1}
  For any $t<\tau\in [0,1]$, we have
  $$
  \hat{D}_tX_\tau(x)=J_\tau K_t\sigma(X_t(x))
  $$
  almost surely, where $ \hat{D}_tX_\tau(x)$ is defined as
  \beaa
(\hna_hX_\tau(x),z)_{\R^n}&=&\int_0^1(\hat{D}_t X_\tau(x),z\otimes \dot{h}(t))_{\R^n\otimes\R^d}dt\\
  &=&\int_0^\tau(\hat{D}_t X_\tau(x),z\otimes \dot{h}(t))_{\R^n\otimes\R^d}dt\,,
  \eeaa
  for $h\in H$, $z\in \R^n$.
  \end{corollary}
\nproof
We can suppose without loss of generality $b=0$, from the relation (\ref{deriv_1}), using Euclidean coordinate
representation we deduce the equation
\begin{equation}
  \label{deriv_2}
  \hat{D}_t^jX^i_\tau=\sigma^i_j(X_t)+\int_t^\tau (\partial\sigma)^i_{kl}(X_s)\hat{D}^j_tX^k_s dm^l_s\,,
\end{equation}
where the sum is performed on repeated indices.
We will show that $(J_\tau K_t\sigma(X_t(x)),\tau\in [t,1])$ satisfies the equation (\ref{deriv_2}). We have
\beaa
\sigma_j^i(X_t)&+&\int_t^\tau(\partial\sigma)_{kl}^i(X_s)[J_\alpha^k(s) K_\beta^\alpha(t)\sigma_j^\beta(X_t)]dm^l_s\\
&=&\sigma_j^i(X_t)+\int_t^\tau(\partial\sigma)_{kl}^i(X_s)[J_\alpha^k(s)dm^l_s K_\beta^\alpha(t)\sigma_j^\beta(X_t)]\\
&=&\sigma_j^i(X_t)+[J_\alpha^i(\tau)-J_\alpha^i(t)] K_\beta^\alpha(t)\sigma_j^\beta(X_t)]\\
&=&\sigma_j^i(X_t)+J_\alpha^i(\tau)K^\alpha_\beta(t)\sigma^\beta_j(X_t)-\delta_{i,\beta}\sigma^\beta_j(X_t)\\
&=&J_\alpha^i(\tau)K^\alpha_\beta(t)\sigma^\beta_j(X_t)\,.
\eeaa
Hence by the uniqueness of the solutions of the equation (\ref{deriv_1}) or (\ref{deriv_2}), the claim follows.
\nqed

\section{\bf{Calculation of Integrands in Martingale Representation and its  Applications }}
\noindent
Here are some applications of the results of the preceding
section. Let us begin with  a non-trivial extension of It\^o-Clark-Hausmann-Bismut-Ocone formula:
\begin{theorem}
\label{Ito-Clark}
Assume that $F(X)\in \M_{2,1}$, then it can be represented as 
$$
F(X)=E[F(X)]+\int_0^1 P(X_s)E[\hat{D}_sF(X)|\calF_s(X)]dB_s\,,
$$
where $\hat{D}_sF(X)$ is defined as
$\hna F(X)(\cdot)=\int_0^{\cdot} \hat{D}_sF(X)ds$ and $P(X_s)$ is the orthogonal projection from $\R^d$ onto $\sigma(X_s)^\star(\R^n)$.
\end{theorem}
\nproof
We know from Theorem \ref{adapt_thm} that $F(X)$ can be represented as 
$$
F(X)=E[F(X)]+\int_0^1P(X_s)\alpha_s(X).dB_s\,,
$$
for some $\alpha(X)\in L^2_a(dt\times d\mu,\R^d)$, adapted to the
filtration $(\calF_t(X),t\in [0,1])$. Moreover, if $F(X)\in \calS(X)$,
then regarding $F(X)$ as a function of the Brownian path $B$, from classical representation theorem (cf. \cite{O1,O2,ASU,ASU-1}), we have 
\begin{equation}
\label{clark}
F(X)=E[F(X)]+\int_0^1 E[D_sF(X)|\calB_s]dB_s\,,
\end{equation}
where $(\calB_s,s\in [0,1])$ is the filtration of the Brownian
motion. It follows from Theorem \ref{iden-thm}, taking the
conditional expectations of both sides of the equality (\ref{clark}),
that
$$
F(X)=E[F(X)]+\int_0^1 P(X_s)E[D_sF(X)|\calF_s(X)]dB_s\,.
$$
Hence in this case $P(X_s)\alpha_s(X)=P(X_s)E[D_sF(X)|\calF_s(X)]=P(X_s)E[E[D_sF(X)|\calF_1(X)|\calF_s(X)]=P(X_s)E[\hat{D}_sF(X)|\calF_s(X)]$
$ds\times d\mu$-a.s. Choose now a sequence $(F_n(X),n\geq 1)\subset
\calS(X)$ which approximates $F(X)$ in $\M_{2,1}$. Then $(\hna
F_n(X),n\geq 1)$ converges to $\hna F(X)$ in $L^2(\mu,H)$ (in other
words $(\hna F_n)$ converges to $\hna F$ in $L^2(\nu,H)$),
consequently 
$$
\lim_nE\int_0^1|P(X_s)E[\hat{D}_sF_n(X)-\hat{D}_sF(X)|\calF_s(X)]|^2ds=0\,,
$$
hence the result follows for any $F(X)\in \M_{2,1}$.
\nqed

\begin{corollary}
\label{ergo_cor}
Assume that $F=F(X)\in \M_{2,1}$ such that $\hna F=0$ $\nu$-a.s. (or
$\hna F(X)=0$ $\mu$-a.s), then $F$ is $\nu$-a.s. (or $F(X)$ is $\mu$-a.s.) a constant.
\end{corollary}
\nproof
The proof follows directly from Theorem \ref{Ito-Clark}.
\nqed

\noindent
We obtain the Poincar\'e inequality from Theorem \ref{Ito-Clark},
whose proof is trivial since the conditional expectation in the
representation formula (\ref{clark}) operates as a contraction:
\begin{theorem}
\label{Poincare}
For any $F(X)\in \M_{2,1}$, we have 
$$
E[|F(X)-E[F(X)]|^2]\leq E[|\rP(X)\hna F(X)|_H^2]\,.
$$
\end{theorem}
\noindent
\begin{theorem}
\label{log_sob}
For any $F\in \M_{2,1}$, it holds true that 
\begin{equation}
\label{log_sob_ineq}
E[F^2(X)\log F^2(X)]\leq E[F^2(X)]\log E[F^2(X)]+2E[|\rP(X)\hna
F(X)|_H^2]\,,
\end{equation}
where $\rP(X)$ is the action induced on $H$ by $(P(X_s),s\in [0,1])$, namely
$$
\rP(X)h(t)=\int_0^t P(X_s)\dot{h}(s)ds\,,
$$
for $h\in H$ and $t\in [0,1]$.
\end{theorem}
\nproof
It suffices to prove the claim for any $F(X)\in \M_{2,1}$ which is strictly
positive and lower bounded by a positive constant. Taking
$f(X)=F^2(X)/E[F^2(X)]$, it suffices to show that 
$$
E[f(X)\log f(X)]\leq \half E\left[\frac{1}{f(X)}|\rP(X)\hna f(X)|_H^2\right]\,.
$$
Using Theorem \ref{Ito-Clark}, we can represent $f(X)$ as
\beaa
f(X) &=&1+\int_0^1 P(X_s)E[\hat{D}_s(f(X))|\calF_s(X)]\cdot dB_s\\
&=&1+\int_0^1
\frac{P(X_s)E[\hat{D}_s(f(X))|\calF_s(X)]}{E[f(X)|\calF_s(X)]}E[f(X)|\calF_s(X)]\cdot
dB_s\,.
\eeaa
Hence $f(X)$ can be represented as an exponential martingale's final
value:
$$
f(X)=\exp\left[\int_0^1\frac{P(X_s)E[\hat{D}_s(f(X))|\calF_s(X)]}{f_s(X)}\cdot
  dB_s-\half
  \int_0^1\left|\frac{P(X_s)E[\hat{D}_s(f(X))|\calF_s(X)]}{f_s(X)}\right|^2 ds\right]\,,
$$
 where $f_s(X)=E[f(X)|\calF_s(X)]$. Let $\ga$ be the probability
 measure on the path space defined by $d\ga=f(X)d\mu$, then the
 Girsanov Theorem and Bayes' formula entail that the relative entropy of $\ga$
 w.r.t. $\mu$, denoted by $H(\ga|\mu)=\int_W \log \frac{d\ga}{d\mu}d\ga$ is equal to
\beaa
H(\ga|\mu)&=&\int \log f(X)d\ga\\
&=&\half
E_\ga\int_0^1\left|\frac{P(X_s)E[\hat{D}_s(f(X))|\calF_s(X)]}{f_s(X)}\right|^2ds\\
&=&\half E_\ga\int_0^1\left|\frac{E[P(X_s)\hat{D}_s(f(X))|\calF_s(X)]}{f_s(X)}\right|^2ds\\
&=&\half E_\ga\int_0^1\left|E_\ga[P(X_s)\hat{D}_s(\log f(X))|\calF_s(X)]\right|^2ds\\
&\leq&\half E_\ga\int_0^1\left|P(X_s)\hat{D}_s(\log f(X))|\right|^2ds\\
&=&\half E\left[\frac{1}{f(X)}|\rP(X)\hna f(X)|_H^2\right]\,.
\eeaa
\nqed

\begin{remark}
Maybe it is pedagogically interesting to note this inequality in terms
of the measure induced by the the diffusion process on the path space,
namely $\nu$:
\begin{equation}
\label{log_sob_2}
\int f^2\log f^2 d\nu-\left(\int f^2d\nu\right)\log \int f^2d\nu\leq 2\int |\rP
\hna f|_H^2d\nu\,
\end{equation}
for any $f\in \M_{2,1}$.
\end{remark}

\section{\bf{Examples and Further Applications}}

\noindent
We shall give applications of the analysis that we have developped  to
some specific cases beginning from the simplest to more
complicated situations

\subsection{Classical Situation}
Let $(X_t,t\geq 0)$ be the solution of the equation
$$
X_t(x)=x+\int_0^t b(X_r(x))dr+\Sigma B_t\,,
$$
where $B=(B_t,t\geq 0)$ is an $d$-dimensional Brownian motion,
$b:\R^n\to \R^n$ is a Lipschitz-continuous vector field and $\Sigma\in
L(\R^n,\R^d)$. 
Regularizing $b$ with a mollifier, we obtain a smooth and Lipschitz
vector field $b_\eps$, let $X^\eps$ be the  solution of SDE
$$
X^\eps_t(x)=x+\int_0^t b_\eps(X^\eps_r(x))dr+\Sigma B_t\,.
$$
For any
$h\in H$, where $H=H^1(\R_+,\R^d)$ is the Cameron-Martin space of $B$, we have
$$
\nabla_hX^\eps_t=\int_0^t\partial b_\eps(X_r)\nabla_hX^\eps_rdr+
\Sigma h(t)
$$
and the Gronwall lemma implies that, for any $T>0$,
$$
\sup_{t\leq T}\|\nabla X^\eps_t\|_{H\otimes\R^n}\leq \|\Sigma\|e^{TK}\,,
$$
where $K$ is the Lipschitz constant of $b$ and $\|\Sigma\|$ is the
operator norm of $\Sigma$. We can pass to the limit
as $\eps\to 0$ to obtain also 
$$
\sup_{t\leq T}\|\nabla X_t\|_{H\otimes\R^m}\leq \|\Sigma\| e^{TK}\,.
$$
The following theorem is now a straight application of Theorem \ref{log_sob}
\begin{theorem}
  \label{log_sob1}
  Let $T>0$ be arbitrary and denote by $\nu^T$  the law of the process
  $(X_t,t\in [0,T])$. Then, for any smooth, cylindrical function on $W_T=C([0,T],\R^d)$, we have
$$
 \int_{W_T}F^2(x)\log\frac{F^2(x)}{\nu^T(f^2)}\nu^T(dx)\leq
 2\|\Sigma\|e^{TK}\int_{W_T}|\nabla F(x)|^2\nu^T(dx)\,,
 $$
in particular, for any $t\in (0,T]$, denoting by $\nu_t$ the law of
the random variable $X_t$, we have
 $$
 \int_{\R^m}f^2(y)\log\frac{f^2(y)}{\nu_t(f^2)}\nu_t(dy)\leq
 2\|\Sigma\|e^{TK}\int_{\R^m}|\partial f(y)|^2\nu_t(dy)\,,
 $$
 for any smooth function $f:\R^n\to \R$.
 \end{theorem}

\section{\bf{Log-Sobolev for Dyson's Brownian Motion}}
 \noindent
 In the example above, the drift coefficient has  nice
 properties and the log-Sobolev inequality is almost trivial. Let us
 look at a more singular case, namely, let $(X_t=(X^1_t,\ldots,X_t^d),\,t\geq 0)$ be the
 Dyson Brownian motion, i.e., the strong solution of the following
 stochastic differential equation:
 \begin{equation}
   \label{Dysoneq}
dX_t^i=dB_t^i+\ga\sum_{1\leq j\neq i\leq
  d}\frac{1}{X_t^i-X^j_t}dt\,,\,i=1,\ldots,d;\,\,\ga>0,
\end{equation}
with $X_0=x_0$ such that $x_1<\ldots<x_d$
We refer the reader to references \cite{CEP,cepa-lep1,cepa-lep2} for
the existence of the strong solutions, where they use the theory of
dissipative mappings in their construction and to \cite{Dy} for the
origins of the equation. For notational convenience, we shall write
above equation as
\begin{equation}
 \label{moneq}
dX_t=\ga m(X_t)dt+ dB_t\,,X_0=x_0\in D(m), 
\end{equation}
where $m$ is the mapping from a subset of $\R^d$, noted as $D(m)$ into
$\R^d$ whose components are given in the equation (\ref{Dysoneq}). In
fact $m$ is a set-valued, dissipative map. It is maximal and this implies that its values reduce to one-point
subsets of $\R^d$ (cf. \cite{Yos}). As a typical example for such a
function, one can take the gradient of a concave function on $\R^d$.
Equation \ref{Dysoneq} corresponds to the case $m=\partial f$, where
$f(x)=\sum_{i<j}\log(x_j-x_i)$ if  $x_1<\ldots<x_d$ and
$f(x)=-\infty$ otherwise. Let us announce first an important technical
result:
\begin{lemma}
  \label{Lip-lemma}
  Let  $X=(X_t,t\geq 0)$ be  the unique strong solution of SDE \ref{moneq}, then, for any $T>0$, $h,k\in
  H([0,T],\R^d)$  we have
  $$
  \sup_{s\leq t}\left|X_s(w+h)-X_s(w+k)\right|\leq 2\sqrt{t}\ga|h-k|_{H([0,t],\R^d)}
  $$
  P-almost surely.
  Consequently, for any $t>0$, $X_t\in\cap_p\DD_{p,1}(\R^d)$, in fact
  $|\nabla X_t|_{\R^d\otimes H([0,t],\R^d)}\in L^\infty(P)$ with the bound
  $2t\ga$, where $H([0,t],\R^d)$ (respectively $H([0,T],\R^d)$ denotes the Cameron-Martin space for
  $\R^d$-valued functions on $[0,t]$ (respectively $[0,T]$).
\end{lemma}
\nproof
We will  suppose that $\ga=1$ and that $(B_t)$ is the linear Brownian
motion, which amounts to working on the canonical Wiener space (for
practical reasons). From the equation \ref{moneq}, it comes
$$
X_t(w+h)=B_t(w)+h(t)+\int_0^tm(X_s(w+h))ds\,,
$$
hence
$$
X_t(w+h)-X_t(w)=h(t)+\int_0^t[m(X_s(w+h))-m(X_s(w))]ds\,.
$$
Let $y_t^h(w)=X_t(w+h)-X_t(w)$, then $(y_t^h,\,t\in [0,T])$ satisfies
$$
y_t^h(w)=h(t)+\int_0^t [m(y^h_s(w)+X_s(w))-m(X_s(w))]ds\,,
$$
hence
\begin{equation}
  \label{fareq}
y_t^h(w)-y_t^k(w)=h(t)-k(t)+\int_0^t[m(y_s^h(w)+X_s(w))-m(y_s^k(w)+X_s(w))]ds\,.
\end{equation}
It follows from equation \ref{fareq} and from the dissipativity of $m$ 
that
\beaa
|y_t^h-y_t^k|^2&=&\int_0^t \frac{d}{ds}|y_s^h-y_s^k|^2ds\\
&=&2\int_0^t(y_s^h-y_s^k,\dot{h}(s)-\dot{k}(s))_{\R^d}+2\int_0^t(y_s^h-y_s^k,m(y^h_s+X_s)-m(y^k_s+X_s))ds\\
&\leq&2\sup_{s\leq t}|y_s^h-y_s^k|\int_0^t|\dot{h}(s)-\dot{k}(s)|ds\\
&\leq&2\sqrt{t}\sup_{s\leq t}|y_s^h-y_s^k|\left(\int_0^t|\dot{h}(s)-\dot{k}(s)|^2ds\right)^{1/2}\\
\eeaa
Since the last line of above inequality is increasing w.r.t. $t\in
\R_+$, the claim follows. The fact that $X_t\in \cap_{p>1}\DD_{p,1}$ follows from \cite{ASU}, Lemma 2 at
  p.72.
\nqed

\noindent
An immediate consequence of Lemma \ref{Lip-lemma} is Fernique's bound:
\begin{corollary}
  \label{expinteg}
  For any $t>0$, there exists some $\eps=\eps(\ga,t)$, such that
  $$
  E\left[\exp\left(\eps\,\sup_{s\leq t}|X_s|^2\right)\right]<\infty
  $$
  
  \end{corollary}

\begin{remark}
Let $T>0$ be any fixed number and let us denote by $\nu$ the law of
the process restricted to the time horizon $[0,T]$, starting from $x_0$. Evidently $\nu$ is
the unique solution of the martingale problem (cf. \cite{S-V}) associated to the
operator $K$ defined as 
$$
Kf(x)=\sum_{i=1}^d\left[\frac{1}{2}\frac{\partial^2f}{\partial
  x_i^2}(x)+\ga m_i(x)\frac{\partial f}{\partial x_i}(x)\right]\,.
$$

It is immediate to see that under $\nu$, the coordinate map $x\to
x(t)$ of $C_{x_0}([0,T],\R^d)$ can be represented as
$$
dx_t=\ga m(x_t)dt+ d\beta_t\,,
$$
where $(\beta_t, t\in [0,T])$ is a $\nu$-Brownian motion. 
\end{remark}

The following is now clear:
\noindent
\begin{proposition}
  \label{closed}
  The covariant derivative defined on the paths of Dyson process
  restricted to $[0,T]$, i.e.,  for $F(X)=f(X_{t_1},\ldots,X_{t_n})$, 
  $$
  \hna_h F(X)=E[\nabla_hF(X)|\calF_T(X)]\,,
  $$
  where $h\in H([0,T],\R^d)$ and  $\hna F(X)=\sum_i\hna_{h_i}h_i$,
 such that  $(h_i)$ is a basis in $H([0,T],\R^d)$, 
  are   closable operators on $L^p(\calF_T(X))$ for any
  $p>1$. Consequently, we have the following representation theorem:
  $$
  F(X)=E[F(X)]+\int_0^TE[D_sF(X)|\calF_s(X)]\cdot dB_s\,.
  $$
  \end{proposition}

From Proposition \ref{closed} we obtain as explained in the general
case, the logarithmic Sobolev inequality:
\begin{proposition}
  \label{Dlog-sob}
  For any $\calF_T(X)$ measurable, cylindrical function $F(X)$, we have
  $$
  E[F^2(X)\log( F^2(X)/E[F^2(X)])]\leq
  2\ga\sqrt{T}E[|\hna F(X)|_{H([0,T],\R^d)}^2]\,.
  $$
\end{proposition}

\section{\bf{Intertwining Relations and Logarithmic Sobolev Inequality}}

\noindent
Let us come back to the general Markovian  situation studied in the paper: if $f$
be a smooth function  on $\R^n$ and if $(Q_t,t\in [0,1])$ is
the semi-group associated to the Markov process induced by $(X_t,t\in
[0,1])$, then $(Q_{1-t}f(X_t(x)),t\in [0,1])$ is an
$(\calF_t(X))$-martingale and, due to the It\^o formula, it can be
represented as
\begin{equation}
  \label{mart1}
  Q_{1-t}f(X_t(x))=Q_1f(x)+\int_0^t\partial_xQ_{1-s}f(X_s(x))\cdot\sigma(X_s(x))dB_s\,.
\end{equation}
From Theorem \ref{Ito-Clark},  we obtain a second representation:
\begin{equation}
  \label{mart2}
  Q_{1-t}f(X_t(x))=Q_1f(x)+\int_0^t P(X_s(x))E[\hat{D}_s
  f(X_1(x))|\calF_s(X)]\cdot dB_s,
\end{equation}
where $\hat{D}_s f(X_1(x))=E[D_sf(X_1(x))|\calF_1(X)]$. Ito Isometry
principle implies that
$$
\partial_xQ_{1-s}f(X_s(x))\cdot\sigma(X_s(x))=
P(X_s(x))E[\hat{D}_sf(X_1(x))|\calF_s(X)]
$$
$ds\times dP$-almost surely, where
$$
\hat{D}_tf(X_1)=E[D_tf(X_1)|\calF_1(X)]\,.
$$
It will be tempting to compare in detail these two representations: As
$$
\hat{D}_tf(X_1)=(\hat{D}_tX_1)^\star\partial f(X_1)\in \R^d\,,
$$
we obtain
\begin{theorem}
  \label{int-thm}
  For any smooth function $f$ on $\R^n$, we have the following
  commutation-intertwining relation:
  \begin{equation}
    \label{int1}
\sigma(X_t)^\star\partial Q_{1-t}f(X_t)=P(X_t) E\left[(\hat{D}_t X_1)^\star\partial f(X_1)|\calF_t(X)\right]\,.
    \end{equation}
\end{theorem}

\noindent
Corollary \ref{calcul1} gives that
$$
(\hat{D}_t X_1)^\star=\sigma^\star(X_t)K_t^\star J_1^\star\,,
$$
where $J_t=\hat{\partial}X_t(x)$ and $ K_t=J_t^{-1}$. Consequently, we
get
\beaa
E[\hat{D}_t
f(X_1)|\calF_t(X)]&=&E[(\hat{D}_tX_1)^\star\hat{\partial}f(X_1)|\calF_t(X)]\\
&=&E[\sigma^\star(X_t)K_t^\star  J_1^\star\hat{\partial}f(X_1)|\calF_t(X)]\\
&=&\sigma^\star(X_t)K_t^\star E[
J_1^\star\hat{\partial}f(X_1)|\calF_t(X)]\,.
\eeaa
\noindent
To calculate the conditional expectation of the last line, we use the
flow property of the diffusion process: For $s<t\leq  1$, let
$X^s_t(x,\omega)$ denote the solution of the SDE defining the
diffusion process at the instant $t$, starting from the point $x\in
\R^n$  at the instant $s$. The flow property implies that
$X_1(x,\omega)=X^t_1(X_t(x,\omega),\tilde{\omega})$, where
$\tilde{\omega}$ denotes a Brownian path independent of
$\{B_\tau:\tau\leq t\}$ Taking the derivative of both sides, from the
chain rule,  it comes
$$
\partial_x X_1(x,\omega)=\partial
X^t_1(X_t(x,\omega),\tilde{\omega})\partial X_t(x,\omega)
$$
almost surely. Hence
\beaa
E[J_1^\star\partial f(X_1)|\calF_t(X)]&=&E[\partial X_1^\star\partial f (X_1)|\calF_t(X)]\\
&=&E[(\partial X^t_1(X_t,\tilde{\omega})\partial X_t(x))^\star\partial
f(X^t_1\circ X_t)|\calF_t(X)]\\
&=&E\left[\partial X_t^\star E[(\partial
X^t_1)^\star\circ(X_t,\tilde{\omega})\partial
f(X^t_1\circ X_t)|\calB_t]|\calF_t(X)\right]\,,
\eeaa
where $(\calB_t,t\geq 0)$ represents the filtration of the underlying
Brownian motion. Due to the independence and the stationarity of the Brownian increments,
the inner conditional expectation can be written as:
\beaa
E[(\partial
X^t_1)^\star\circ(X_t,\tilde{\omega})\partial
f(X^t_1\circ X_t)|\calB_t]&=& E[(\partial
X^t_1)^\star(\xi)\partial
f(X^t_1(\xi))]|_{\xi=X_t(x)}\\
&=&E[(\partial
X_{1-t})^\star(\xi)\partial
f(X_{1-t}\xi))]|_{\xi=X_t(x)}
\eeaa
and in particular, it is $\calF_t(X)$-measurable. Let us define the
semigroup $(M_t,t\geq 0)$ on the vector-valued functions as
$$
M_t\alpha(\xi)=E[(\partial X_t(\xi))^\star\alpha(X_t(\xi))]\,,
$$
$\xi\in \R^n$. Using this semigroup we end up with the identity
\bea
\label{meq1}
 E[J_1^\star\partial f(X_1)|\calF_t(X)]&=&E[\partial X_t^\star|\calF_t(X)]\,M_{1-t}\partial
f(X_t(x))\\
&=&\hat{\partial} X^\star_t(x) M_{1-t}\partial f(X_t(x))\,.\nonumber
\eea
Finally we find that
\begin{equation}
  \label{meq2}
E[\hat{D}_t f(X_1(x))|\calF_t(X)]=\sigma(X_t(x))^\star M_{1-t}\partial f(X_t(x))
\end{equation}
almost surely. Although the expression (\ref{meq2}) is what we can
obtain by differentiation under the expectation, it has an important
byproduct that we write separately as a lemma:
\begin{lemma}
  \label{mart_lemma}
  For any smooth function $f$ on $\R^n$, the process $(\hat{\partial}
  X^\star_t(x) M_{1-t}\partial f(X_t(x)),\,t\in [0,1])$ is an $\calF_\cdot(X(x))$-martingale.
\end{lemma}

The following result may be called as a modified Poincar\'e
inequality:
\begin{theorem}
  \label{mod_poincare}
  Let $Z$ be the functional defined as
  $$
  Z(x)=\int_0^1|\sigma^\star(X_t(x))(\hat{\partial}
  X_t(x)^\star)^{-1}|^2dt\,.
  $$
  We then have, for any $f\in C_b^2(\R^n)$, the following inequality:
  \begin{equation}
    \label{poin_ineq}
    E[|f(X_1(x))-E[f(X_1(x))]|^2]\leq
    E\left[Z\left(\hat{\partial}X_1(x)\cdot\hat{\partial}X_1^\star(x)\partial
    f(X_1(x)),\partial f(X_1(x))\right)_{\R^n}\right]\,.
    \end{equation}
    \end{theorem}
    \begin{remark}
    Let operator-valued function $\Gamma(x,y)$ be defined as
    \begin{equation}
      \label{Gamma}
    \Gamma (x,y)=E[Z \hat{\partial} X_1(x)\cdot\hat{\partial} X_1^\star(x)|X_1(x)=y]\,,
    \end{equation}
    then we can read the same inequality in $\R^n$ as
    $$
    \int_{\R^n}\left|f(y)-\int_{\R^n}f(z)q_1(x,z)dz\right|^2q_1(x,y)dy\leq
    \int_{\R^n}\left(\Gamma (x,y)\partial f(y),\partial f(y)\right)_{\R^n}q_1(x,y)dy\,,
    $$
    where $y\to q_1(x,y)$ is the density of the law of $X_1(x)$.
  \end{remark}
  
    \nproof
    From the martingale representation theorem, from the relation
    (\ref{meq2}) and from the It\^o isometry, we obtain
    \beaa
 E[|f(X_1(x))-E[f(X_1(x))]|^2]&=&E\int_0^1|\sigma^\star(X_t(x))
 M_{1-t}\partial f(X_t(x))|^2dt\\
 &=&E\int_0^1|\sigma^\star(X_t(x)) (\hat{\partial} X^\star_t)^{-1}(\hat{\partial} X^\star_t)
 M_{1-t}\partial f(X_t(x))|^2dt\\
 &\leq&E\int_0^1|\sigma^\star(X_t(x)) (\hat{\partial} X^\star_t)^{-1}|^2|(\hat{\partial} X^\star_t)
 M_{1-t}\partial f(X_t(x))|^2dt\\
 &\leq&E\int_0^1|\sigma^\star(X_t(x)) (\hat{\partial} X^\star_t)^{-1}|^2|(\hat{\partial} X^\star_1)
 \partial f(X_1(x))|^2dt\,,
 \eeaa
where the inequality at the last line  follows from Lemma
\ref{mart_lemma}.
\nqed

\noindent
Same observation gives also an extension of the logarithmic Sobolev
inequality for the measure $q_1(x,y)dy$  which is the law of the random variable
$\omega\to X_1(\omega, x)$:
\begin{theorem}
  \label{log_sob_1}
  Let $f\in C^2_b(\R^n)$, then, we have
  \begin{equation}
    \label{ineq4}
    E\left [f^2(X_1)\log \frac{f^2(X_1)}{E[f^2(X_1)]} \right]\leq
    2\int_{\R^n}\left(\Gamma(x,y) \partial f(y), \partial f(y)\right)_{\R^n}
    q_1(x,y)dy\,,
  \end{equation}
  where $\Gamma$ is defined with (\ref{Gamma}).
\end{theorem}
\nproof
We proceed as in the proof of Theorem \ref{log_sob} with $f\in
C_b^2(\R^n)$ which is strictly positive. Moreover we can suppose that
$E[f(X_1(x))]=1$. In this case the claimed inequality takes the form
$$
 E\left [f(X_1)\log f(X_1)\right]\leq
    \int_{\R^n}(\Gamma(x,y) \partial f(y), \partial f(y))_{\R^n}\frac{1}{f(y)}
    q_1(x,y)dy\,.
 $$
Let $d\ga=f(X_1(x))dP$, similar calculations of entropy $H(\ga|P)$ as in the proof of Theorem
\ref{log_sob} imply that
\beaa
H(\ga|P)&=&\int\log f(X_1(x)) d\ga\\
&=&\half E_\ga\int_0^1\left|\frac{P(X_s)E[\hat{D}_sf(X_1)|\calF_s(X)]}{Q_{1-s}f(X_s)}\right|^2\\
&=&\half
E_\ga\int_0^1\left|\frac{P(X_s)\sigma^\star(X_s(x))M_{1-s}\partial
    f(X_s(x))}{Q_{1-s}f(X_s)}\right|^2\\
&\leq&E_\ga\int_0^1\left|\frac{\sigma^\star
    (X_s) (\hat{\partial}X_s^\star)^{-1}\hat{\partial} X_s^\star
    M_{1-s}\partial f(X_s)}{Q_{1-s}f(X_s)}\right|^2ds\\
&\leq&E_\ga\int_0^1|\sigma^\star(X_s)
(\hat{\partial}X_s^\star)^{-1}|^2 \left|\frac{\hat{\partial} X_s^\star
    M_{1-s}\partial f(X_s)}{Q_{1-s}f(X_s)}\right|^2ds\,.
\eeaa
Note that, from Lemma \ref{mart_lemma}, the process
$$
(\frac{\hat{\partial} X_s^\star M_{1-s}\partial
  f(X_s)}{Q_{1-s}f(X_s)}, s\in [0,1])
$$
is a $\ga$-martingale, consequently we have
\beaa
H(\ga|P)&\leq&E_\ga\left[\left|\frac{\hat{\partial} X_1^\star\partial
      f(X_1)}{f(X_1)}\right|^2\int_0^1|\sigma^\star(X_t)
  (\hat{\partial}X_s^\star)^{-1}|^2dt\right]\\
  &=&E\left[\frac{1}{f(X_1)}|\hat{\partial}X_1^\star\partial f(X_1)|^2Z\right]\,,
  \eeaa
  where $Z$ is defined in the announcement of Theorem \ref{mod_poincare}.
  \nqed

\noindent
Left side of the equation (\ref{ineq4}) depends only on the law of the
random variable $X_1(x)$, but at the right hand side there is the
term $\Gamma(x,y)$ which is defined as  
$$
 \Gamma (x,y)=E[Z \hat{\partial} X_1(x)\cdot\hat{\partial}
 X_1^\star(x)|X_1(x)=y]\,,
 $$
 where
 $$
 Z(x)=\int_0^1|\sigma^\star(X_t(x))(\hat{\partial}
  X_t(x)^\star)^{-1}|^2dt\,.
  $$
Consequently, the right hand side of the inequality (\ref{ineq4}) depends rather
explicitly on $\sigma$, while the left hand side depends only on
$a=\sigma\sigma^\star$, hence only on the law of the
process. Meanwhile, this inequality remains valid for any $\sigma$
which is $n\times d$-matrix valued function satisfying the relation
$a=\sigma\sigma^\star$. In particular the second dimension $d$ is not
fixed. Let us define
$$
\Sigma_a=\cup_{d\geq 1}\{\sigma \in C_b^1(\R^n, M(n\times
  d)):\,\sigma\,\sigma^\star=a\}\,,
  $$
where $M(n\times d)$ denotes $n\times d$-matrices.
  \begin{corollary}
    \label{final}
    Let us define the functional $Q(\partial f,\partial f)$ for $f\in
    C_b^2(\R^n)$ as
    $$
    \inf_{\sigma\in \Sigma_a}\int_{\R^n} (\Gamma_\sigma (x,y)\partial
    f(y),\partial f(y))_{\R^n}q_1(x,y)dy\,.
    $$
   Let $\nu_1^x$ be the law of the random variable $\omega\to X_1(x,\omega)$. Then it holds true that
    \begin{equation}
      \label{ineq5}
      \int_{\R^n}f^2(y)\log\frac{f^2(y)}{\nu_1^x(f^2)}\nu_1^x(dy)\leq
      Q(\partial f,\partial f)\,.
    \end{equation}
  \end{corollary}

  \begin{remark}
    Assume that, for any smooth, upper and lower bounded function $f$  on
    $\R^n$, we have the following property: the process
    \begin{equation}
      \label{ls_cond}
  (t,w)\to \frac{|\sigma(X_t)^\star\partial
    Q_{1-t}f(X_t)|^2}{Q_{1-t}f(X_t)}
  \end{equation}
  is a submartingale, then, replacing $f$ by $\sqrt{f}$,  we have the following version of
  logarithmic Sobolev inequality, which follows from the It\^o formula:
  $$
  E\left[f(X_1)^2\log \frac{f(X_1)^2}{E[f(X_1)]^2}\right]\leq 2E
  [|\sigma^\star(X_1)\partial f(X_1)|^2]\,,
$$
for any smooth function $f$.
\end{remark}

  \section{\bf{Conditionally Orthogonal Functionals and Inequalities}}
  \noindent
  Let us begin with the general case of the stochastic differential equation without Markov
  coefficients, i.e., the $b$ and $\sigma$ may depend in a causal manner on the Wiener path
  explicitly.

  Let us define by $R_s(X)$ the orthogonal projection defined as
  $$
  R_s(X)=I_{\R^d}-P_s(X)\,,
  $$
  and let $\rR(X)$ be its action on the $H$-valued functions:
  $$
  \rR(X) h(t)=\int_0^tR_s(X)\dot{h}(s)ds\,.
  $$
    
  This projection operator gives birth to a new concept of stochastic process which is intimately
  related to the degeneracy of the diffusion process $X=(X_t,t\in [0,1])$. Let us begin by
  \begin{theorem}
    \label{ortho_mart}
    Let $u\in L^2_a(\mu,H)$, define
    $$
    L=\exp\left(\int_0^1(R(X_s)\dot{u}_s,dB_s)-\half\int_0^1|R(X_s)\dot{u}_s|^2ds\right)
    $$
    and assume that $E_\mu[L]=1$. Then we have
    \begin{equation}
      \label{c_Girsanov}
      E_\mu[L|\calF_1(X)]=1
    \end{equation}
    $\mu$-almost surely. In particular, if $u$ is $\calF_1(X)$-measureable, then 
    $$
    E\left[\exp\left(\int_0^1(R(X_s)\dot{u}_s,dB_s)\right)|\calF_1(X)\right]=
    \exp\half\int_0^1|R(X_s)\dot{u}_s|^2ds,
    $$
    hence under the conditional probability $\mu^X$, defined by $\mu^X(\cdot)=\mu(\cdot|\calF_1(X))$,
    $$
    \int_0^1(R(X_s)\dot{u}_s,dB_s)
    $$
    is a Gaussian random variable with mean zero and variance $|\rR(X)u|_H^2$. In particular, the process
    $$
    N_t=\int_0^tR_s(X)dB_s
    $$
    is an additive (i.e., with independent increments) Gaussian martingale under the conditional probability $\mu^X$.
  \end{theorem}
  \nproof
  Let $L_t=\exp\left(\int_0^t(R(X_s)\dot{u}_s,dB_s)-\half\int_0^t|R(X_s)\dot{u}_s|^2ds\right)$, from the It\^o formula, we have
  $$
  L_t=1+\int_0^t L_s(R_s(X)\dot{u}(s),dB_s).
  $$
From Theorem \ref{iden-thm}, we  have 
\begin{equation}
  \label{ec}
E[L_1|\calF_1(X)]=1+\int_0^1 (P_s(X),R_s(X)E[\dot{u}(s)L_s|\calF_s(X)],dB_s)=1
\end{equation}
almost surely. If $\dot{u}$ is adapted to the filtration $\calF(X)$, from the relation (\ref{ec}), we obtain
$$
E\left[\exp\la\int_0^t(R(X_s)\dot{u}_s,dB_s)|\calF_t(X)\right]=\exp\frac{\la^2}{2}\int_0^t|R_s(X)\dot{u}_s|^2ds
$$
a.s., for any $\la\in \R$. The last claim follows from the fact that, for any $h,k\in H$,
$$
E[\exp\int_0^1(R_s(X)\dot{h}(s)+R_s(X)\dot{k}(s),dB_s)|\calF_1(X)]=
\exp\half \int_0^1(|R_s(X)\dot{h}(s)|^2+|R_s(X)\dot{k}(s)|^2)ds
$$
a.s. if $\dot{h}$ and $\dot{k}$ have disjoint supports.
\nqed

\begin{proposition}
  \label{cond_prop}
  Let $F=F(X,N)$ be in $L^2(\mu)$ be a cylindrical function, define $\nabla^{(2)}_hF(X,N)$ as
  $$
  \nabla^{(2)}_hF(X,N)=\lim_{\eps\to 0}\frac{F(X,N+\eps h)-F(X,N)}{\eps}\,,
  $$
  for $h\in H$ and $D_\tau^{(2)}F(X,N)$ as
$$
  \nabla^{(2)}_hF(X,N)=\int_0^1(D_\tau^{(2)}F(X,N),\dot{h}(\tau))_{\R^d}d\tau.
$$
We have then
\begin{equation}
  \label{itrt}
F(X,N)=E^X[F]+\int_0^1(R_\tau(X)E^X[D^{(2)}_\tau F(X,\cdot)|\mathcal{N}_s],dB_s)\,,
\end{equation}
where $(\mathcal{N}_s,s\in [0,1])$ is the filtration generated by $N=(N_s,s\in [0,1])$ and $E^X$ is
the conditional expectation corresponding to $\mu^X$. Consequently the following conditional
log-Sobolev inequality holds true:
$$
E^X\left[F^2\log\frac{F^2}{E^X[F^2]}\right]\leq 2 E^X[|\rR\nabla^{(2)}F|_H^2]\,.
$$
\end{proposition}
\nproof
As the process $N$ is conditionally an additive Gaussian process it has martingale representation property under the law $\mu^X$. Besides, from Theorem \ref{c_Girsanov}, we obtain the following conditional integration by parts formula 
\begin{equation}\label{ibp_1}
E^X\left[F(X,N)\int_0^1(R_s(X)\dot{u}_s(X,N),dB_s)\right]=
E^X\left[\int_0^1(R_s(X)\dot{u}_s(X,N),D_s^{(2)}F(X,N))ds\right]\,,
\end{equation}
for any bounded, and $\mathcal{N}$-adapted process $\dot{u}(X,\cdot)$. This implies the representation (\ref{itrt} by conditioning. 
The rest of  proof follows from the relation (\ref{ibp_1}) by using the usual techniques.
\nqed

{\bf\noindent{Acknowledgment:}} This work has been done in the  University of Bilkent in Ankara. Moreover its actual  form is due to   a very careful reading of an  anonymous referee to whom the author is grateful.

\vspace{2cm}
\footnotesize{
\noindent
A. S. \"Ust\"unel, Bilkent University, Math. Dept., Ankara, Turkey\\
ustunel@fen.bilkent.edu.tr}

\end{document}